\documentclass[12pt,openbib]{article}
\usepackage{amsmath}
\usepackage[cp1251]{inputenc}
\usepackage[english]{babel}
\usepackage{amsfonts}
\usepackage{amsfonts,amssymb}
\usepackage{amssymb}
\usepackage{latexsym}
\usepackage{euscript}
\usepackage{enumerate}
\usepackage{graphics}
\usepackage[dvips]{graphicx}

\newtheorem{lem}{\bf Lemma}
\newtheorem{rem}{\bf Remark}
\newtheorem{prop}{\bf Proposition}

\newtheorem{df}{\bf Definition}

\begin{document}

\begin{center}
 {\large\bf On Long Virtual Biquandles}  \vspace{12mm}

D. A. Fedoseev \vspace{3mm}

\textit{Moscow State University, Main Building, \\ Chair of Differential Geometry and Applications, \\ Moscow, 119991, Leninskie Gory, 1, Russia} \vspace{5mm}

\textbf{Abstract} 

Virtual quandles with two operations are discussed in the article. Certain knot invariant is constructed and used to distinguish two long virtual knots. \vspace{3mm}


\end{center}

\vspace{3mm}

\section{Introduction}

Object named "quandle" is well-known in modern knot theory. It provides good knot invariants. We will remind how this object can be constructed (as described in [1]).

Let $\Gamma$ be a finite set if "colours" with an operation "circle": $\circ \colon \Gamma \times \Gamma  \to \Gamma.$

\begin{df}
\textit{Correct colouring} of an oriented knot (link) diagram $D$ is such correlation between arcs of the diagram $D$ and elements of $\Gamma,$ that for each crossing the following is verified: $c = a \circ b$ if arcs are marked with colours $a, b$ and $c$ as shown on the diag.1:
\end{df}

\begin{center}
\includegraphics[width=30mm]{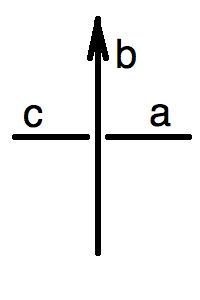}
\end{center}

\begin{rem}
We do not look at the orientation of arcs denoted a and c.
\end{rem}

Now we will enforce several conditions on the operation "circle" which ensure that the number of correct colourings of an oriented diagram is invariant under Reidemeister movements. Direct computation shows that the conditions are as follows:

\begin{enumerate}
\item $\forall a \in \Gamma$ $a \circ a = a;$ 
\item $\forall a,b \in \Gamma$ equation $x \circ a = b$ has exactly one solution $x \in \Gamma.$ Further it will be denoted as $b/a;$ 
\item $\forall a,b,c \in \Gamma$ $(a \circ b) \circ c = (a \circ c) \circ (b \circ c).$
\end{enumerate} 

Any set with operation "circle" satisfying the above conditions is called a \textit{quandle}.

From the very definition follows

\begin{prop}
The number of correct colourings of an oriented diagram with elements of a quandle is a knot (link) invariant.
\end{prop}

There also exists a more general approach to the construction of a quandle, the one using generators and relations.

Let $A$ be an alphabet \--- a set of letters. A \textit{word} in alphabet $A$ is by definition any finite sequence of elements of $A$ and symbols $\circ$ and $/$. Now we will define a set $D(A)$ of \textit{allowed words}. $D(A)$ is defined inductively, following the following rules:

\begin{enumerate}
\item Any letter of the alphabet $A$ is an allowed word;
\item If $W_1, W_2 \in D(A)$, then $(W_1) \circ (W_2)$ and $(W_1) / (W_2)$ are allowed words;
\item There are no other allowed words.
\end{enumerate}

Further throughout the text we will ignore writing brackets in cases when the meaning of the structure is clear.

Consider a set of relations $R=\{r_\alpha = s_\alpha|r_\alpha, s_\alpha \in D(A)\}.$ We introduce an equivalence relation for $D(A)$ such that for any $W_1, W_2 \in D(A)$ $W_1 \equiv W_2$ if ad only if there exists a sequence of transformations beginning with $W_1$ and ending $W_2$ constructed according to the following rules (trivial equivalences):

\begin{enumerate}
\item $x \circ x \Leftrightarrow x;$
\item $(x \circ y)/y \Leftrightarrow x;$
\item $(x/y) \circ x \Leftrightarrow x;$
\item $(x \circ y) \circ z \Leftrightarrow (x \circ z) \circ (y \circ z);$
\item $r_i \Leftrightarrow s_i.$
\end{enumerate}

A set of allowed words factorized according to this equivalence is, clearly, a quandle with operation $\circ.$

Now for a given knot we construct a quandle invariant according to the following scheme. First of all wi assign a letter to each arc of the knot diagram and take this set of letters as an alphabet. Then we produce a set of relations $R:$ for every crossing of the diagram we state $a \circ b = c$ (as shown on the diag.1). After that we construct a quandle as described above. Such "\textit{knot quandle}" is an almost complete knot invariant in the sense that if two knots are equal, corresponding quandles are isomorphic. It is not very convenient one, though, because it is usually difficult to verify if two quandles are isomorphic or not. So some modifications of the structure are considered ang used.

\section {Basic constructios}

Let's consider a virtual knot. An object not unlike quandle can be constructed \--- a \textit{virtual quandle}.

For now we will consider "long arcs" of a virtual knot diagram \--- a connected component of the set obtained from the diagram by deleting all virtual crossings. Again we label all the long arcs with letters (generators) $x_i$ and note the same relations $a \circ b = c$ for each classical crossing of the diagram (again we consider an oriented knot or link).

\begin{df}
A quandle (according to Kauffman) \--- is a formal quandle of a knot, obtained by ignoring all the virtual crossings of the diagram.
\end{df}

The object defined above provides some knot invariants but it is comparatively weak. For example virtual trefoils (right and left) cannot be distinguished with it.

A better generalization of a quandle is called a \textit{virtual quandle}.

\begin{df}
A virtual quandle is a quandle $(M, \circ)$ with an operation $f$ such that there exists inverse operation $f^{-1}$ and for any $a,b \in M$ $f(a) \circ f(b) = f(a \circ b).$
\end{df}

Now we construct an invariant $Q(L)$ for a given oriented diagram $L$ of a virtual knot $K$. The structure present will be a strong virtual knot invariant.

First of all we label all the arcs of the diagram with letters $a_i.$ Let $X(L)$ be a set of words obtained inductively using letters $a_i$ and symbols $\circ,$ $/$, $f$, $f^{-1}.$ We will factorize this set using a following equivalence: a transitive and reflexive closure of the following set of trivial equivalences:

\begin{enumerate}
\item $f^{-1}(f(a)) \Leftrightarrow f(f^{-1}(a)) \Leftrightarrow a;$
\item $x \circ x \Leftrightarrow x;$
\item $(x \circ y)/y \Leftrightarrow x;$
\item $(x/y) \circ x \Leftrightarrow x;$
\item $(x \circ y) \circ z \Leftrightarrow (x \circ z) \circ (y \circ z);$
\item $f(a \circ b) \Leftrightarrow f(a) \circ f(b);$
\end{enumerate}

furthermore for every classical crossing we state $a_{i_1} \circ a_{i_2} \Leftrightarrow a_{i_3},$ as shown on the diag.2:

\begin{center}
\includegraphics[width=30mm]{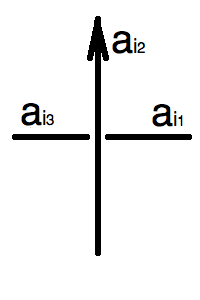}
\end{center}

and for every virtual crossing we state $x' \Leftrightarrow f(x)$ and $y' \Leftrightarrow f(y)$ as shown on the diag.3:

\begin{center}
\includegraphics[width=30mm]{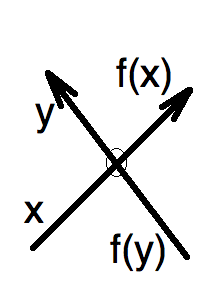}
\end{center}

\begin{prop}
Virtual quandle $Q(L)$ constructed as show above is an invariant of virtual knots (links).
\end{prop}

Rigor proof is given in [2].

Until now all the classical crossings were "equal" in the sense that we applied the same equivalence for arcs incident to any classical crossings. If we can somehow divide all the classical crossings into two categories, we can construct a \textit{biquandle} $(M, \circ, \star)$ which gives stronger invariant than the one described above. Good example of such an object can be presented using long virtual oriented knots.

\begin{df}
Long virtual biquandle is a set $M$ with operations $\circ,$ $\star,$ $\bar{\circ},$ $\bar{\star}$ and $f$ such that $(M, \circ, f)$ is a virtual quandle, $(M, \star, f)$ is a virtual quandle and the following is verified:

$\cdot$ $\forall a,b\in M$ $(a\circ b)/_{\circ}b=(a/_{\circ}b)\circ b=(a\star b)/_{\star}b=(a/_{\star}b)\star b=a$,

$\cdot$ $\forall a,b,c\in M$  $ (a\diamond b)\bullet c=(a\bullet c)\diamond(b\bullet c)$, where $\diamond$ and $\bullet$ \--- are some operations from the following list: $\circ$, $\star$, $/_{\circ}$, $/_{\star}$,

$\cdot$ $\forall a,b\in M$ $ f(a\diamond b)=f(a)\diamond f(b),$ where $\diamond$ \--- is an operation from the list $\circ$, $\star$,

$\cdot$ $\forall a\in M$ $ f(f^{-1}(a))=f^{-1}(f(a)) = a.$

And "strange relations" are verified:

$\cdot$ $\forall x,a,b\in M$ $ x\diamond(a\circ b)=x\diamond(a\star b)$,

$\cdot$ $\forall x,a,b\in M$ $ x\diamond(a/_{\circ}b)=x\diamond(a/_{\star}b)$,

where again $\diamond$ are some operations from the list $\circ$, $\star$, $/_{\circ}$, $/_{\star}$.
\end{df}

Now for a given diagram we produce a free long biquandle (a quandle, formally generated by operations $\circ,$ $\star,$ $\bar{\circ},$ $\bar{\star}$ and factorized according to relations $1$ \--- $5$) and then we factorize it according to the structure of the diagram: we state a relation $c = a \circ b$ for every classical crossing which is an \textit{early overcrossing} according to the knot's orientation; $c = a \star b$ for every early undercrossing and treat virtual crossings as above (as shown on the diag.3).

This object gives us a knot invariant.

\section{Construction of a long virtual biquandle}

To construct an example of a long virtual biquandle invariant we will use the following fact. Let $G$ be a group such that $\exists a, b \in G:$ $[a,b] \notin Z(G)$ but there exists $n \in \mathbb{N}$ such that for any $a,b \in G$ $[a,b^n] \in Z(G)$ but $\exists a, b \in G:$ $[a, b^n] \neq e.$ Here square brackets denote a commutator in the group (i.e. $[a,b]=ab-ba$), $Z(G)$ is the group's center and $e$ denotes neutral element in the group. If given such a group we can use $G$ as an alphabet and define $a \circ b := bab^{-1}$ and $a \star b := b^{n+1}ab^{-n-1}.$ Operation $f$ can be chosen freely (though it must "respect" both binary operations $\circ$ and $\star$).

Let us give an example of such a group using a Cayley graph. The graph we will be using is a square divided into $64$ smaller squares; horizontal sides of those a marked with letter $a$, vertical \--- with letter $b.$ Finally all the horizontals and verticals of the big square are oriented: the lowest horizontal is oriented right, next one left and so on; the leftmost vertical is oriented up, next one \--- down and so on.

Here $a$ and $b$ are generators and relations are given by the graph, assuming that the square is glued into a torus and all the horizontals an verticals are oriented as described.

The group $G$ consists of all "paths" in the graph. To elements are considered equal if the corresponding paths connect the same vertices.

\begin{lem}
a) $\exists x, y \in G: [x,y] \notin Z(G);$

b) $\forall x, y \in G [x,y^2] \in Z(G)$ and $\exists x, y \in G: xy^2x^{-1}y^{-2} \neq e.$
\end{lem}
$\square$
a) Let $x=a,$ $y=b.$

Then $a(aba^{-1}b^{-1}) = a^3b^2$ and $(aba^{-1}b^{-1})a = a^3b^{-2} = a^3b^6 \neq a^3b^2.$ Therefore $[a,b] \notin Z(G).$

b) Obviously, every $x$ in $G$ is equal to an element of the form $a^kb^l.$ So we are to prove that for any $i, j, k, l$

$$A=(a^{k}b^{l})(a^{i}b^{j})^{2}(a^{k}b^{l})^{-1}(a^{i}b^{j})^{-2}=a^{k\pm i\pm i\pm(-k)\pm(-i)\pm(-i)}b^{l\pm j\pm j\pm(-l)\pm(-j)\pm(-j)} = a^\alpha b^\beta\in Z(G)$$
and $\exists i, j, k, l:$ $A \neq e.$

$A$ depends solely on parity of numbers $i, j, k, l.$ Direct computation shows the following correspondence between parity of $i, j, k, l$ and numbers $\alpha$ and $\beta:$

\begin{center}
\begin{tabular}{|c||c|c|c|c|c|c|c|c|c|c|c|c|c|c|c|c|}
\hline
$i$ & $1$ & $1$ & $0$ & $0$ & $1$ & $1$ & $0$ & $0$ & $1$ & $1$ & $0$ & $0$ & $1$ & $1$ & $0$ & $0$ \\
\hline
$j$ & $1$ & $0$ & $1$ & $0$ & $1$ & $0$ & $1$ & $0$ & $1$ & $0$ & $1$ & $0$ & $1$ & $0$ & $1$ & $0$ \\
\hline
$k$ & $1$ & $1$ & $1$ & $1$ & $1$ & $1$ & $1$ & $1$ & $0$ & $0$ & $0$ & $0$ & $0$ & $0$ & $0$ & $0$ \\
\hline
$l$ & $1$ & $1$ & $1$ & $1$ & $0$ & $0$ & $0$ & $0$ & $1$ & $1$ & $1$ & $1$ & $0$ & $0$ & $0$ & $0$ \\
\hline
\hline
$\alpha$ & $0$ & $-4i$ & $0$ & $0$ & $0$ & $0$ & $0$ & $0$ & $0$ & $-4i$ & $0$ & $0$ & $0$ & $0$ & $0$ & $0$ \\
\hline
$\beta$ & $0$ & $0$ & $-4j$ & $0$ & $0$ & $0$ & $-4j$ & $0$ & $0$ & $0$ & $0$ & $0$ & $0$ & $0$ & $0$ & $0$ \\
\hline
\end{tabular}
\end{center}

Therefore not every possible value of $A$ is equal to $e,$ but certainly all the values of $A$ are in $Z(G).$ So the lemma is proved.
$\blacksquare$

Finally, we construct the biquandle $(G, \circ, \star, f):$ $$x \circ y = yxy^{-1};$$ $$x \star y = y^3xy^{-3};$$ $$f(a) = ab; f(b) = b;$$ $$\forall \alpha \in G \quad f(\alpha^{-1})=f(\alpha)^{-1};$$ $$\forall \alpha, \beta \in G \quad f(\alpha \beta) = f(\alpha)f(\beta).$$

\section{Use of the quandle}

For example we will use the biquandle constructed above to distinguish right and left long virtual trefoils. We will use "colourings invariant" with elements of the biquandle. It is important to notice that the number of correct colourings with the colour of the first (according to orientation) long arc of a long knot fixed is invariant. Even more, in that case the set of possible colourings of the second long arc is invariant as well.

$\includegraphics[width=70mm]{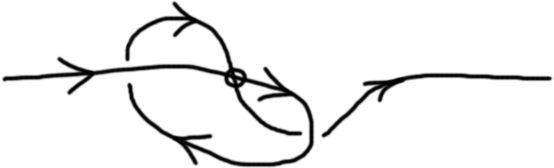}$
$\includegraphics[width=70mm]{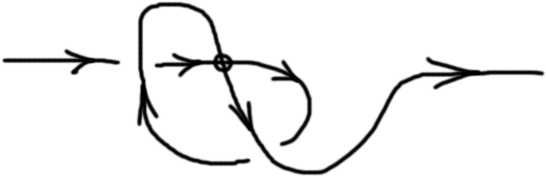}$

Let arcs of the first knot be labeled $a_i$ (according to orientation) and arcs of the second knot be labeled $b_i.$ Let $a_1 = a.$

In that case we have: $a_2 = ab^{-1}, a_3 = a^2b^{-1}a^{-1}, a_4 = (ab)^2a^{-1}, a_5 = ab^2.$

To show inequality of the knots it is enough to prove that there is no correct colouring of the second knot with $b_1 = a, b_5 = ab^2.$

Assume that is not the case:  $b_1 = a, b_5 = ab^2.$ Then $b_4 = ab, b_2 \star b_4 = a, b_4 \star b_5 = b_3, b_2 = f(b_3).$ Therefore $\alpha = (ab)^{-3}a(ab)^3=(ab^3)^3ab(ab^3)^{-3} = \beta.$ But direct computions show taht $\alpha = a^7b^6$ and $\beta = a^7b.$ So $\alpha \neq \beta$ and our assumption is incorrect. So we have proved the inequality of the knots under consideration.

This work is supported by RFFI grant (project ~10-01-00748-a), grant of President of RF: aid for Scientific Schools (project 3224.2010.1),
program Development of Scientific Potential of Higher School (project
2.1.1.3704), programs Scientific and Scientifically-Teaching personnel of innovate Russia (contracts 02.740.11.5213 and 14.740.11.0794).

\end{document}